\documentclass[a4paper]{amsart}  

\usepackage{amssymb} \usepackage{amscd} \usepackage{times}

\def\zbb{\mathbb{Z}}  
  
  \def\phi{\varphi}
 \def\p1{{\mathbb{P}^1_\zbb}}

\newtheorem{Theorem}{\quad Theorem}[section]

\newcommand{\be} {\begin{equation}}
\newcommand{\ee} {\end{equation}}

\begin{document}

\title{ An inequality for the solutions to an equation with boundary condition.}

\author{Samy Skander Bahoura}

\address{Department of Mathematics, Pierre and Marie Curie University, 75005, Paris, France.}
              
\email{samybahoura@gmail.com} 

\date{}

\maketitle

\begin{abstract}

We give a uniform estimate and an inequality for solutions of an equation with Dirichlet boundary condition.

\end{abstract}

\section{Introduction and Main Result}

We set $ \Delta = -\nabla^i \nabla_i$  the Laplace-Beltrami operator on a connected compact Riemannian manifold with boundary $ (M,g)$  of dimension $ n \geq 3 $ with metric $ g $.

\smallskip

We consider the following boundary value problem:

\begin{displaymath}  (P)  \left \{ \begin {split} 
      \Delta u + h u  & = V u^{N-1}, \,\,u >0     \,\, &&\text{in} \!\!&& M, \\
                  u  & = 0  \,\,             && \text{on} \!\!&&\partial M.               
\end {split}\right.
\end{displaymath}

Here, we assume the solutions in the sense of distributions and here also in $ C^{2,\alpha}(\bar M), \alpha >0 $ and the operator $ \Delta +h $ is coercive in $ H^1_0(M) $ with $ h $ a smooth function and $ 0 \leq V \leq b $ with $ V\not \equiv 0 $  a smooth function and $ ||V||_{C^{\alpha}} \leq A $ and $ N=\frac{2n}{n-2}$ the critical Sobolev exponent. This problem arises in physics and astronomy.

\smallskip

In  [1-25] we have various  estimates and inequalities of type $ \sup +\inf$ and $ \sup\times \inf $. Here we look to the case of lower bounds for theses inequalities as in [2], [4], [5], [7] and [21].

\smallskip

Here we prove an a priori estimate and an implicit Harnack type inequality for the solutions of the problem $(P)$.

Our main results are:

\begin{Theorem} For all $ b, A,  M_0 >0 $ and all compact $ K $ of $ M $, there is a positive constant $ c=c(b, A, \alpha, K, M, g, n, M_0)>0 $ such that:

$$ \inf_K u \geq c \,\, {\rm if} \,\, \sup_M u\leq M_0. $$

for all solution $ u >0 $ of $ (P)$ relative to $ V\geq 0 $ with bounds $ b, A $.
\end{Theorem}

In the previous Theorem, we have an apriori estimate for the solutions of the problem $ (P) $. We can see, the influence of the boundary condition, the Dirichlet boundary condition.

\smallskip

As a consequence of the previous Theorem, we have $ \inf_K u \geq c=c(\sup_M u) $ with $ c>0 $. We write this as: the following implicit Harnack type inequality:

\begin{Theorem} For all $ b, A >0 $ and all compact $ K $ of $ M $, there is a positive constant $ c>0 $ depending on the $ \sup_M u $, such that:

$$ \inf_K u \geq c(b, A, \alpha, K, M, g, n, \sup_M u). $$

for all solution $ u >0 $ of $ (P) $ relative to $ V \geq 0 $ with bounds $ b, A $.
\end{Theorem}

\smallskip

In the last Theorem, we have an implicit inequality between $ \sup $ and $ \inf$: if we have an information on the sup then we have information of the inf and also if we have an idea on the infimum, then we have an information of the supremum.

\smallskip

$$ \sup_M u \,\,{\rm and} \,\, \inf_K u\,\, {\rm are \,\, linked \,\, implicitly}. $$

\smallskip

In various previous papers, we proved some inequalities of type $ \sup \times \inf $ for equations of type prescribed scalar curvature and under some conditions, on the manifold or on the solutions.

\smallskip

In [2, 5], to have a positive lower bound for $ \sup \times \inf $, we also supposed the manifold compact without boundary. We used the Green function of an invertible operator, lower bound of the Green function and an iterate scheme.

\smallskip

In [5], for manifolds with boundary, we have added a condition on the solutions(a jump condition), to prove a minoration of $ \sup \times \inf $. Here we can see the influence of the Dirichlet boundary condtion. We also used, Green function for an invertible operator and used a positive lower bound of the Green function, locally, and an iterate scheme.

\smallskip

In [4], we can see that the optimal $ \sup \times \inf $ inquality holds. Here, we have an implicit Harnack inequality for more general nonlinear equations with critical Sobolev exponent. Here also, we used, Green function and lower bound of the Green function, locally, and two methods to obtain the estimate, one by iterate scheme, the other by integration by parts in the equation, and by the Sobolev inequality.

\smallskip

In [7], we used a blow-up technique to have an explicit Harnack inequality on the $n-$ball of the euclidean space.

\smallskip

{\bf Remarks:}

\smallskip

1) The function $ M_0 \mapsto c $ is non-increasing function of $ M_0 $.

\smallskip

2) We need uniform $ C^{\alpha} $ regularity for $ V $, because we will use the strong maximum principle which require $ C^2 $ regularity of the solutions.

\smallskip

3) We can replace the exponent $ N $ by $ q \geq 2+\epsilon_0 $ with $ \epsilon_0 >0$.

\smallskip

4) These results can be applied  to $ h=\frac{n-2}{4(n-1)}S_g$ with $ S_g $ the scalar curvature, in the non-negative case. We obtain the prescribed scalar curvature equation and the Yamabe equation on manifolds with boundary, with Dirichlet boundary condtions. For example $ S_g \geq 0 $ and $ S_g >0 $ somewhere.

\smallskip

5) In [4], we have explicit Harnack inequality, for $ K $ compact subset of $ \Omega $, there is $ c=c(n,K,\Omega) >0 $ and $ \tilde c=\tilde c(n, K,\Omega) >0 $ such that, precisely:

$$ \sup_{\Omega} u \times \inf_K u \geq c(n, K,\Omega) \int_{\Omega} |\nabla u|^2 dx \geq \tilde c(n, K, \Omega) >0. $$

Where $ u >0 $ is solution of

\begin{displaymath}  \left \{ \begin {split} 
      \Delta u  & = u^{N-1-\epsilon}, \,\,u >0     \,\, &&\text{in} \!\!&& \Omega \subset \subset {\mathbb R}^n, \\
                  u  & = 0  \,\,             && \text{on} \!\!&&\partial \Omega.               
\end {split}\right.
\end{displaymath}

Here $ \epsilon >0 $ small and perhaps tending to $ 0 $. It is a nonlinear boundary Dirichlet problem with subcritical exponent perhaps tending to the critical exponent.

\smallskip

For the previous inequality we used that the solutions are uniformly bounded in a neighborhood of the boundary.(By the moving-planes method).

\smallskip

One can see, in the previous inequality that if the solutions have an upper bound, then locally they have a positive lower bound. 

\smallskip

\section{Proofs of the Theorems}

\smallskip

Proofs of Theorems 1 and 2:

\smallskip

We use a cover argument (a finite cover) for a local convergence in open sets of charts, to have convergence of subsequence and to build a limit function $ u\geq 0 $ by extensions. 

\smallskip

We argue by contradiction, suppose that for a $ K $ compact of $ M $ we have  sequences of solutions $ (u_i,V_i) $ of $ (P) $ such that: $ \inf_K u_i \to 0 $ and $ \sup_M u_i \leq M_0 $. Using charts around points of $ \bar M $ in  ${\mathbb R}^n $ and $ {\mathbb R}^n_+ $, we can reduce the uniform estimates to estimates in the euclidean case and by the elliptic estimates, see [17] , we have $ u_i\to u \geq 0 $ uniformly in $ C^{2,\alpha}(\bar M) $. (Since $ u_i $ is regular, we have for a chart of a boundary, $ \psi $, $ u_io\psi $ is regular in a half ball, and the usual trace can be considered and we reduce the problem to a boundary value problem with a portion boundary as trace's subset, then, theorems of [17] can be applied to $ u_io\psi $, chapters 9 and 6).

\smallskip

One can build a limit function by using local convergence in open sets of charts. Since, all  the topologies coincide, the topology of charts and the topology of the metric.

\smallskip

Now we use the condition $ \inf_K u_i \to 0 $ to have for the solution $ u \geq 0 $, $ u(x_0)=0 $ with $ x_0 \in K \subset \subset M $.

\smallskip

Then we use  the maximum principle applied to $ -u $ to have $ u\equiv 0 $ on $ \bar M $. (since the operator$ \Delta +h $ is coercive, we have the existence of the Green function, $ G $, of $ \Delta+h$ with Dirichlet boundary conditions. We use the Green representation formula to prove $ u\equiv 0 $, because $ Vu^{N-1}\equiv 0 $ and $ G >0 $ in $ M $, and we multiply the equation of $ u $ by $ u $ and we integrate, and then, we use coercivity to have $ u\equiv 0 $). 

Also, we can use the fact $ u\geq 0$ we take $ \lambda <0$ such that $ -h+\lambda <0$, the operator $ -\Delta -h+\lambda $ satisfy the strong maximum principle and we have $ (-\Delta -h+\lambda)(-u)=(Vu^{N-2}-\lambda)u \geq 0 $ and $ -u\leq 0$, and $ (-u)(x_0)=0, x_0\in K \subset \subset M $ and then $ -u\equiv 0$, then $ u\equiv 0$.

Thus,

$$ u_i \to 0 \,\,{\rm in }\,\, C^{2,\alpha}(\bar M). $$

But if we multiply the equation by $ u_i $ and use integration by part, the coercivity of the operator and the Sobolev embedding we obtain for some constant $ C_2>0 $ independant from $ i $:

$$ C_2||u_i||_N^2 \leq C_1||u_i||_{H_0^1(M)}^2 \leq \int_M |\nabla u_i|^2+hu_i^2  $$

and,

$$\int_M |\nabla u_i|^2+hu_i^2 = \int_M (\Delta u_i + h u_i)u_i  = \int_M V_i u_i^N \leq b ||u_i||_N^N, $$

thus,

$$ ||u_i||_N^{N-2} \geq \frac{C_2}{b}>0, $$

with uniform bound. 

\smallskip

And this contradict: $ u_i\to 0 $ in $ C^0(\bar M) $.

\smallskip

Thus, for all $ b, A >0 $, $ M_0 >0 $, for all compact $ K $ of $ M $, there is a positive constant $ c=c(b, A, \alpha, K, M, g, n, M_0) $ such that:

$$ \sup_M u \leq M_0 \Rightarrow  \inf_K u \geq c >0, $$

with by definition, 

$$ c=c(M_0)=\sup \{ c'>0, \inf_K u \geq c', \,{\rm if} \, \sup_M u\leq M_0, M_0 >0 \}. $$

With this definition, by taking twice the supremum in $ c'>0 $, we see that if $ 0 < \sup_M u \leq M_0\leq K_0 \Rightarrow \inf_K u \geq c(K_0) \Rightarrow \inf_K u\geq c(M_0) \geq c(K_0) \Rightarrow M_0\mapsto c $ is non-increasing function of $ M_0 $.

\smallskip

Take $ u >0 $ a regular solution of $ (P) $ relative to $ V $ and set $ M_0= \sup_M u >0 $, then for all $ v >0 $ solution of $(P)$ relative to $ W\geq 0 $ with bounds $ b, A $,  we have the previous inequality.

\smallskip

But $ v=u >0 $ is a particular solution, thus it satisfies the inequality and we have:

$$ \inf_K u \geq c=c(b, A, \alpha, K, M, g, n, \sup_M u)>0.$$

Thus, we obtain an implicit Harnack type inequality.

\bigskip

\underbar { \bf Questions:} 1) How to expilcit the previous Harnack type inequality ? in the previous papers we proved that we have inequalities of type $ \sup \times \inf $.

\smallskip

2) Can we have inequalities of type:

$$ (\sup_M u)^{\alpha} \times \inf_K u \geq c >0, \alpha >0 ? $$

\end{document}